\documentclass{amsart}
\usepackage[numbers]{natbib}
\usepackage{geometry}
\usepackage{fancyhdr}
\usepackage{afterpage}
\usepackage{longtable}
\usepackage{lscape}
\usepackage{array}
\usepackage{color}
\usepackage{tabularx}
\usepackage{ltablex}
\usepackage{graphicx}
\usepackage{amsmath,amssymb,amsbsy}
\usepackage{amsfonts}
\usepackage{amssymb}
\usepackage{amsthm}
\usepackage{dcolumn,array}
\usepackage{algpseudocode}

\usepackage[pageanchor=true,plainpages=false,pdfpagelabels,bookmarks,bookmarksnumbered]{hyperref}
\usepackage{enumerate}
\usepackage{rotating}
\usepackage[all]{xy}

\newcommand{\g}{{\text{gal}}}

\newcommand{\Z}{{\mathbb Z}}

\newcommand{\Q}{{\mathbb Q}}

\newcommand{\Aut}{\text{Aut}}
\DeclareMathOperator{\Gal}{Gal}

\usepackage[T1]{fontenc}

\begin{document}

\title{Finding Galois splitting models to compute local invariants}
\author{Benjamin Carrillo}
\maketitle

\begin{abstract}
For prime $p$ and small $n$, Jones and Roberts have developed a database recording invariants for $p$-adic extensions of degree $n$. We contributed to this database by computing the Galois slope content, Galois mean slope, and inertia subgroup for a variety of wildly ramified extensions of composite degree using the idea of \emph{Galois splitting models}. We will describe a number of strategies to find Galois splitting models including an original technique using generic polynomials and Panayi's root finding algorithm.
\end{abstract}

\section{Introduction}

For a number field $K$, i.e. a finite extension of $\Q$, and for each prime $p$, we have an associated $p$-adic algebra $K \otimes \Q_p \cong \prod^g_{i=1} K_{p,i}$, where each $K_{p,i}$ is a finite extension of $\Q_p$. \\

We can answer a variety of questions about $K$ using basic invariants of the $K_{p,i}$, such as their ramification index and residue field degree. For more advanced questions, there is a need for more detailed information about the $p$-adic extension, such as their Galois group and ramification groups, which allows us to measure the wild ramification of the extension.\\
 
Let $\mathcal{K}(p, n)$ be the set of isomorphism classes of degree $n$ extensions of $\Q_p$ which is finite, as a consequence of Krasner's Lemma \cite{lang}. We want to record and store information that encapsulates the filtration of the Galois group by its ramification groups for each extension for particular $\mathcal{K}(p, n)$, as well as some companion results. To easily compute this data we use the idea of \emph{Galois splitting models}. \\

Formally, a \textit{Galois splitting model} of a $p$-adic extension $K/\Q_p$, is a polynomial $f(x)\in \Z[x]$ that is irreducible in $\Z_p[x]$ with $\Q_p(\alpha) \cong K$ and $\Gal(\Q_p(\alpha)/\Q_p) =\Gal(\Q(\alpha)/\Q)$ where $\alpha$ is a root of $f(x)$. With a Galois splitting model $f(x)$ the computation of various invariants related to $K$ that we are interested in can be computed easily using $f(x)$ and its corresponding number field.\\

For a given $p$-adic extension we use techniques coming from inverse Galois theory such as class field theory and generic polynomials to find a Galois splitting model. In particular, we will describe an original technique to find our desired polynomial using generic polynomials and Panayi's root finding algorithm. \cite{pauli}\\

In keeping with Jones and Roberts \cite{jj_local}, all data will be available at \url{https://hobbes.la.asu.edu/LocalFields/} as well as \url{http://www.lmfdb.org/LocalNumberField/} so all computations are recorded once and is freely available for those who are interested. Also an implementation of the main algorithm described here will be available at \url{https://github.com/bcarril1/gsm_panayi}.\\ 

\section{Preliminaries}
\subsection{Panayi's Root Finding Algorithm}
We will now describe Panayi's root finding algorithm, let $\varphi(x) \in \mathcal{O}_K[x]$, where $K$ is a finite $p$-adic extension with uniformizer $\pi$. \cite{pauli} \cite{panayi}.\\

\begin{itemize}
\item Let $\varphi(x)=c_nx^n+ c_{n-1}x^{n-1}+ \ldots +c_1x+c_0 \in \mathcal{O}_K[x]$ and define the valuation $\nu_K(\varphi):=\min\{\nu_K(c_n), \ldots, \nu_K(c_0)  \}$ where the initial $\nu_K$ is the $\pi$-adic valuation of $K$ and $\varphi^{\#}(x):= \varphi(x)/\pi^{\nu_K(\varphi)}$. For $\alpha \in \mathcal{O}_K$ denote its representative in the residue field $k$ by $\overline{\alpha}$, and for $\beta\in k$, denote a lift of $\beta $ to $\mathcal{O}_K $ by $\hat{\beta}$.\\

\item To count the number of root of $\varphi(x)$ in $K$ we define two sequences $(\varphi_i(x))_i$ and $(\delta_i)_i$.\\

\item Set $\varphi_0(x):= \varphi^{\#}(x)$ and let $\delta_0 \in \mathcal{O}_K$ be the lift of a root of $\overline{\varphi_0(x)}$.\\

\item If $\overline{\varphi_i^{\#}(x)}$ has a root $\beta_i$ then define $\varphi_{i+1}(x):= \varphi_i^{\#}(x\pi+\hat{\beta_i})$ with $\delta_{i+1}:= \hat{\beta_i}\pi^{i+1}+\delta_i$.
\end{itemize}
At some point, one of the following cases must occur:
\begin{enumerate}
\item $\deg(\overline{\varphi_i^{\#}})=1$ then $\delta_{i-1}$ is an approximation of one root of $\varphi(x)$
\item $\deg(\overline{\varphi_i^{\#}})=0$ then $\delta_{i-1}$ is not an approximation of a root of $\varphi(x)$
\item $\overline{\varphi_i^{\#}}$ has no roots and thus $\delta_{i-1}$ is not an approximation of a root of $\varphi(x)$
\end{enumerate}
If at any step of the process there exist multiple roots $\beta_i$ for $\overline{\varphi_i(x)}$, we split the sequence and proceed.\\

\subsection{Parametric/Generic polynomials} The Inverse Problem of Galois Theory asks if for a finite group $G$ and a field $K$, does there exist a finite Galois extension $L$ such that $G=\Gal(L/K)$. An additional question is if $G$ can be realized as the Galois group of a field extension of $K$, can we construct a family of polynomials over $K$ such that the Galois group of the polynomials over $K$ is $G$.\\

The idea of \emph{parametric polynomials} is an attempt to answer the second question. Following Jensen, Ledet, and Yui, consider the polynomial  $P(\mathbf{t},x) \in K(\mathbf{t})[x]$, where $\mathbf{t}=(t_1, \ldots, t_n)$ and the $t_i$ are indeterminants. Let $\mathbb{L}$ be the splitting field of $P(\mathbf{t},x)$ over $K(\mathbf{t})$. We say that $P(\mathbf{t},x)$ parametrizes $G$-extensions and is a \emph{parametric polynomial} if $P$ satisfies the following two properties:
\begin{enumerate}[1.]
\item $\mathbb{L}/K(\mathbf{t})$ is Galois with Galois group $G$  \\
\item Every Galois extension $L/K$ with Galois group $G$ is the splitting field of a polynomial $P(\mathbf{a},x)$ for some $\mathbf{a} \in K^n$.
\end{enumerate}
If $P(\mathbf{t},x)$ has the additional property of parametrizing $G$-extensions for any field containing $K$ then we say $P(\mathbf{t},x)$ is a \emph{generic polynomial}. Discussion on generic polynomials up to degree 7 is readily available, for example, see \cite{jensen2002generic}.\\

For a Galois group $G$ it can be identified as a transitive subgroup of $S_n$. When referring to Galois groups we will use standard notation ($S_n$, $A_n$, $C_n$, $D_n$) as well as T-numbering that was introduced in \cite{butler}, writing $n$T$j$ for a degree $n$ field whose Galois closure has Galois group T$j$ of $S_n$.\\

\subsection{Confirming Galois splitting models}

Given a degree $n$ extension $F/\Q_p$ with defining polynomial $f(x)$ and Galois group $G$ with T-number T$j$, as we produce polynomials we naturally need to confirm if any of those polynomials are a Galois splitting model for $F/\Q_p$. If we have a degree $n$ polynomial $g(x)$ and $K/\Q = \Q[x]/\langle g(x) \rangle$ has Galois group with T-number T$j$, then $g(x)$ is a Galois splitting model for $F/\Q_p$ if $f(x)$ has a root in $\hat{K}/\Q_p$, the completion of $K/\Q$ by the use of Panayi's Algorithm.\\

If $f(x)$ does indeed have a root in $\hat{K}/\Q_p$, this means that $F^\g$ is contained in $\hat{K}^\g$, thus $F^\g \subseteq \hat{K}^\g$. But the T-number of the Galois group of $F/\Q_p$ is T$j$ and the T-number of the Galois group of $\hat{K}/\Q_p$ is at most T$j$, therefore both Galois groups must have the same T-number and $F^\g = \hat{K}^\g$ and thus $F=\hat{K}$ and the T-number of the Galois group of $K$ and $\hat{K}$ are the same.\\

\section{Constructing Galois splitting models}
In this section, we describe four strategies to find candidates for Galois splitting models.

\subsection{Using a Database Search} Our initial attempt to find Galois splitting models is to use the various databases of number fields, namely the databases of Jones and Roberts \cite{jj_global}, Kl{\"u}ners and Malle \cite{kluners}, and the LMFDB ($L$-function and Modular Forms Database) \cite{lmfdb} to find any initial matches. We quickly filter out number fields where the prime $p$ splits within the number field and with the remaining polynomials we check if any are Galois splitting models for some $p$-adic extension. This strategy is useful for finding quick matches in our initial step to find Galois splitting models for $p$-adic extensions of a given degree.

\subsection{Using Galois Theory}
The next strategy for finding Galois splitting models is to use group theoretic facts about the Galois group of a $p$-adic extension to construct a Galois splitting model using composita of smaller fields.\\

Given a degree $n$ field extension $K/F$ with Galois group $G$ we want to determine if there exist a subfield $L$ of $K^\g$ such that $L^\g=K^\g$ and $L$ is the compositum of two smaller subfields. This can be easily found using group theoretic arguments. Namely, we are searching for two non-trivial subgroups of the Galois group of index less than $n$ such that their intersection is trivial. In the case of multiple pairs of non-trivial subgroups that satisfies the previous statement, we pick the pair that generates the largest group, as this will correspond to a common subfield for the fixed fields of the pair of subgroups and we want to minimize this degree.\\

We will show an example of this process using a field extension $K/F$ with Galois group 14T9. Using Magma we find that there exist a pair of subgroups $H$ and $K$ (of index 2 and 8 respectively) whose intersection is trivial and $\langle H,K \rangle$ is the group 14T9. Their corresponding fixed fields will have Galois group 2T1 and 8T25.\\

This means there exists a degree 16 extension of $K^\g$ that is the compostium of a degree 2 and degree 8 extension  with no common subfield. The Galois closure of this degree 16 extension has its Galois group isomorphic to a 14T9 group. We can then use Magma to recover the 14T9 extension using resolvents.\\

Once we identify a suitable pair of subgroups for a particular Galois group for an extension of $\Q_p$, we need to find a Galois splitting model for the fixed field of each subgroup. To find these Galois splitting models with common subfield, we again refer back to the Galois splitting models already found, the various number field databases, or use the other methods which will be described below. We choose to find Galois splitting models for the lower degree extensions as it could potentially easier due to a greater number of lower degree number fields in the various databases and could be quicker  computationally with the methods that will be described next.\\

\subsection{Using Class Field Theory}
This strategy uses class field theory to find Galois splitting models. We know that solvable Galois extensions can be constructed by a chain of abelian extensions and we use this idea to construct a Galois splitting model. \\

For a number field $K$ we can use class field theory to construct a cyclic extension of $K$ with prime order and conductor who divides an ideal of $\Z$ that we specify. We use two implementations one in Pari/GP \cite{PARI2} implementing algorithms from \cite{cohen2012advanced} and the other in Magma \cite{fieker2001computing}.\\

Here is an example with a field extension $K/\Q$ with Galois group 15T26. We cannot find a Galois splitting model for this type of extension using composita of proper subfields of $K^\g$. It can be calculated that $|\Aut_{\Q}(K)|=3$. Thus using standard Galois Theory there exist a unique degree 5 subfield $F$ of $K$, and it can be shown that $\Gal(F/\Q)\cong C_5$. Now $K$ over $F$ is a Galois extension of degree 3, therefore $\Gal(K/F) \cong C_3$. So finding a $C_3$ extension of $F$, will give us a degree 15 extension over $\Q$ that may include or be a 15T26 extension. Thus for a 15T26 $p$-adic extension, we will find a Galois splitting model for the degree 5 subfield  and use this process to try to construct a Galois splitting model for the full extension.\\

\begin{figure}[h]
\[
\xymatrix@R-5pt{
 K\rlap{}\ar@{-}[dd]^{C_3} \ar@{}@/_1.5pc/[dddd]_{\text{15T26}}\\
 \\
F \\
\\
\Q  \ar@{-}[uu]_{C_5}}\]
\caption{Example using 15T26}
\end{figure}
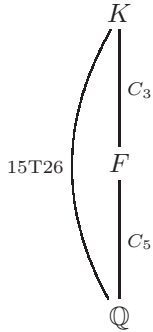

Another more complicated example is a 15T33 extension of $\Q$. Now $|\Aut_{\Q}(K)|=1$ and there exist a degree 5 subfield $F$ with Galois group $C_5$ over $\Q$. So $K/F$ is a degree 3 extension that is not Galois, this means the Galois group of $K/F$ must be $S_3$. The Galois closure $K'$ of $K$ over $F$ has degree 6,  so there exists a degree 2 extension $L_1$ over $F$, which is always cyclic. Additionally, since $K'$ over $F$ is Galois, then $K'$ over $L_1$ is a degree 3 Galois extension, which means it must also be cyclic. \\
 
We can use our process in two steps. The first step is the find quadratic extensions of $F$ and then find $C_3$ extensions of the resulting fields. This will give us degree 30 extensions that could possibly contain a 15T33 extension. From a degree 30 extension we can use Pari/GP or Magma to try to find a 15T33 field extension. A way to optimize this process is to identify that $\Gal(L_1/\Q) \cong C_{10}$, and then find $C_3$ extensions of $L_1$ in one iteration rather than two iterations. Again for a 15T33 $p$-adic extension we will find a Galois splitting model for the $C_{10}$-extension and use this process to try to find a Galois splitting model for the 15T33 $p$-adic extension.\\

\begin{figure}[h]
\[
\xymatrix@R-5pt{
 & K'\rlap{}\ar@{-}[dl]\ar@{-}[ddr]^{C_3} \\
K \\
& & L_1 \\
&F  \ar@{-}[ur]_{C_2} \ar@{-}[uul]\\
& \Q \ar@{-}[u]^{C_5}}\]
\caption{Example using 15T33}
\end{figure}
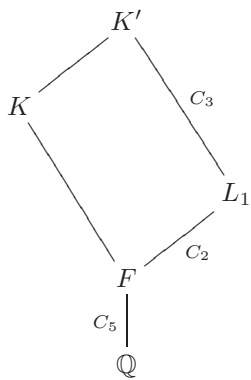

Using group theory we can tell if a number field with particular Galois group $G$ is the subfield of an extension whose Galois group is the wreath product of two groups with the first group cyclic of prime order, i.e. it provides the base extension and order of cyclic extension from the base extension. Kl{\"u}ner and Malle's database provides us this information so we have a starting point to apply this technique to find a Galois splitting model.\\

\subsection{Using Generic Polynomials}
We will describe our strategy to find Galois splitting models using generic polynomials. Given an extension $K/\Q_p$, there will exist a subfield $F$ (which could be trivial) such that $F^{unram}=K^{unram}$, where $K^{unram}$ is the maximal unramified extension of $K$. If there exists multiple subfields of $K$, we simply pick the subfield $F$ such that $[K:F]$ is minimal. If we know $\Gal(K^\g/\Q_p)$, then we can determine $G'=\Gal(K^\g/F)$. Let $P(\mathbf{t},x)$ be generic polynomial that parametrises $G'$-extensions and $f(x)$ be a Galois splitting model for $F/\Q_p$ with $L=\Q[x]/\langle f(x) \rangle$.\\ 

If we calculate $P(\mathbf{a},x)$ for $\mathbf{a} \in L^n$. This will give us a relative field extension over $L$ which we can then view as an extension over $\Q$. Once we calculate a defining polynomial for the absolute extension, we determine if it is a Galois splitting model for full extension $K/\Q_p$.\\

To find suitable values of $\mathbf{a}\in L^n$ when we have an explicit parametric or generic polynomial, we developed an algorithm that uses Panayi's root finding algorithm. The idea of this algorithm is to "zero in" to correct values of the indeterminants of $P(\mathbf{t},x)$ by using Panayi's algorithm and the appropriate substitutions of the indeterminants of $P(\mathbf{t},x)$. We will describe the algorithm now:\\

Let $\pi_F$ be a uniformizer of $ F/\Q_p \cong \hat{L}/\Q_p $ and $\pi_K$ be a uniformizer of $K/\Q_p$. Note we can represent $\pi_F$ in terms of $\pi_K$ by factoring the minimal polynomial of $\pi_F$ over $\Q_p(\pi_K)$.\\

\subsection*{Algorithm 1} 
\begin{itemize}
\item
Let $\mathbf{b}=(b_1,\ldots,b_{n-1},t)$ with $b_i \in L$ and define $$\varphi(x,t):=P(\mathbf{b},x) =c_nx^n+ c_{n-1}x^{n-1}+ \ldots +c_1x+c_0 \in \left(\mathcal{O}_L[t]\right)[x]$$ where $c_i=a_{i,0}+a_{i,1}t+\ldots+a_{i,m_i}t^{m_i}$.

\item
Define $\nu_K(\varphi):=\min\{\nu_K(a_{0,0}), \ldots, \nu_K(a_{n,m_n})  \}$ and $\varphi^{\#}(x,t):= \varphi(x,t)/\pi_K^{\nu_K(\varphi)}$. For $\alpha \in \mathcal{O}_K$ denote its representative in the residue field $k$ by $\overline{\alpha}$, and for $\beta\in k$, denote a lift of $\beta $ to $\mathcal{O}_K $ by $\hat{\beta}$.

\item
We initially create a set $S=\{ \{\{\}, \varphi^{\#}(x,t) \} \}$ and for $\{s,\varphi(x,t)\} \in S$, if at any point $\nu_K(a_{i,0})\geq \min\{\nu_K(a_{i,1}),\ldots , \nu_K(a_{i,m_i})\}$ for some $i$, we then replace $\{s,\varphi(x)\}$ with $\bigcup_{\beta \in k} \{\{s \cup \{\hat{\beta}\}, \, \varphi(x,\hat{\beta}+\pi_F t) \}\}$ in $S$.\\

We choose to make the substitution when $\nu_K(a_{i,0})\geq \min\{\nu_K(a_{i,1}),\ldots , \nu_K(a_{i,m_i})\}$ for some $i$, because otherwise the choice of $t$ may affect the value of $\nu_K(\varphi)$. Thus we substitute all possible values for the first $\pi_F$-adic digit of $t$ and then proceed with the algorithm.\\

\item
If $\overline{\varphi}(x,t)$ has a root $\beta$ then define $\varphi'(x,t):= \varphi^{\#}(x\pi_K+\hat{\beta},t)$ and we replace $\{s,\varphi(x,t)\}$ with $\{ s, \, \varphi'(x,t) \}$.
\end{itemize}

When one of the following cases occur:
\begin{enumerate}
\item $\deg_t(\overline{\varphi}(x,t))=0$ and $\deg_x(\overline{\varphi})=1$ then $P(\mathbf{a},x)$ with $\mathbf{a}=(b_1,\ldots,b_{n-1},\sum_{i=1} s_i\pi_F^i)$ has a root of in $K$.
\item $\deg_t(\overline{\varphi}(x,t))=0$ and $\deg_x(\overline{\varphi})=0$ then $P(\mathbf{a},x)$ with $\mathbf{a}=(b_1,\ldots,b_{n-1},\sum_{i=1} s_i\pi_F^i)$ does not have root of in $K$.
\item $\deg_t(\overline{\varphi}(x,t))=0$ and $\overline{\varphi}(x,t)$ has no roots then $P(\mathbf{a},x)$ with $\mathbf{a}=(b_1,\ldots,b_{n-1},\sum_{i=1} s_i\pi_F^i)$ does not have root of in $K$.
\item
The cardinality of $s$ is larger than a predefined bound.
\end{enumerate}

$\newline$\\
Once a list of polynomials $P(\mathbf{a},x)$ is produced from Algorithm 1, we can then create a list of polynomials over $\Q$. From this list of polynomials we search for a Galois splitting model for the extension $K/\Q_p$. Note we chose to give values from $L$ to all but one indeterminant, but this is not necessary as one can easily modify the algorithm to solve for multiple indeterminants. In general only having one indeterminant does make the computations quicker.\\

\subsubsection{Examples}
For our first example we want to find a Galois splitting model for every $D_5$-extension of $\Q_5$ which there are 3 such extensions. So let $F=\Q_5$ and we choose the Galois splitting model for $F$ to be $f(x)=x$ and hence $\pi_F=5$. From \citep{jensen2002generic} a generic polynomial for $D_5$-extensions is $P(s, t, x) = x^5 + (t - 3)x^4 + (s - t + 3)x^3+(t^2 - t - 2s - 1)x^2 + sx + t$. Since the polynomial $P$ has two parameters, we will choose $s=5$. \\

We will walk through a full example now. For the first extension $K/\Q_5$ its defining polynomial is $g(x) = x^5 + 15x^2 + 5$ and for the first iteration we let $\varphi(x,t) = x^5 + (t - 3)x^4 + (-t + 8)x^3 + (t^2 - t - 11)x^2 + 5x + t$ and $s=\{\}$.\\

Since $\nu_L(a_{0,0})=\nu_L(0)>\nu_L(a_{0,1})=\nu_L(1)$ we will do the substitution step and substitute $3+\pi_F t$ for $t$ and therefore we get $\varphi(x,t) = x^5 + 5tx^4 + (-5t + 5)x^3 + (25t^2 + 25t - 5)x^2 + 5x + (5t + 3)$ and $s=\{3\}$.\\

Now $\overline{\varphi}(x,t)$ has a root of 2, and thus 
$\varphi(x\pi_K+2,t) = (-15\pi_K^2 - 5)x^5 + (5t + 10)\pi_K^4x^4 + (35t + 45)\pi_K^3x^3 + (25t^2 + 115t + 105)\pi_K^2x^2 + (100t^2 + 200t + 125)\pi_Kx + (100t^2 + 145t + 65)$.\\

But $\nu_L(a_{0,0})=\nu_L(65) \geq \nu_L(a_{0,1})=\nu_L(145)$, so for this substitution step we will substitute $\pi_F t$ and we will get $\varphi(x,t) = (-15\pi_K^2 - 5)x^5 + (25t + 10)\pi_K^4x^4 + (175t + 45)\pi_K^3x^3 + (625t^2 + 575t + 105)\pi_K^2x^2 + (2500t^2 + 1000t + 125)\pi_Kx + (2500t^2 + 725t + 65)$
therefore $\varphi'(x,t) = \varphi^{\#}(x\pi_K + 2,t) =  (-3\pi_K^2 - 1)x^5 + (5t + 2)\pi_K^4x^4 + (35t + 9)\pi_K^3x^3 + (125t^2 + 115t + 21)\pi_K^2x^2 + (500t^2 + 200t + 25)\pi_Kx + (500t^2 + 145t + 13)$  and $s=\{3,0\}$.\\

Now $\varphi'(x,t)$ has a root of 3, therefore $\varphi'(x\pi_K+3,t) =(225\pi_K^4 + 150\pi_K^2 + 25)x^5 + (-75\pi_K^4 + (-125t + 3325)\pi_K^3 + (5625t + 2250)\pi_K^2 + 1125\pi_K + (1875t + 750))x^4 + ((-4500t - 1800)\pi_K^4 + (-2625t - 1125)\pi_K^3 + (-1500t + 19650)\pi_K^2 + (-875t - 225)\pi_K + 6750)x^3 + ((625t^2 + 575t - 3945)\pi_K^4 + (-20250t - 8100)\pi_K^3 + (-23625t - 7425)\pi_K^2 + (-6750t - 2700)\pi_K + (-7875t - 2025))x^2 + ((4725t + 1215)\pi_K^4 + (3750t^2 + 3450t - 5445)\pi_K^3 + (2500t^2 - 39500t - 16075)\pi_K^2 - 2025\pi_K + (-13500t - 5400))x + ((2025t + 810)\pi_K^4 + (4725t + 1215)\pi_K^3 + (5625t^2 + 5175t - 2700)\pi_K^2 + (7500t^2 + 3000t + 375)\pi_K + (2500t^2 + 725t - 1150))$.\\

And $\varphi''(x,t) = \varphi'^{\#}(x\pi_K + 3,t)= (-3\pi_K^4 - \pi_K^2)x^5 + (-45\pi_K^3 + (-75t - 30)\pi_K^2 - 15\pi_K + (-25t - 10))x^4 + ((60t + 24)\pi_K^4 + (35t + 9)\pi_K^3 - 270\pi_K^2 - 90)x^3 + (54\pi_K^4 + (270t + 108)\pi_K^3 + (315t + 81)\pi_K^2 + (125t^2 + 115t + 21)\pi_K)x^2 + ((-100t^2 - 40t - 5)\pi_K^4 + 81\pi_K^3 + (540t + 216)\pi_K^2 + (-1500t^2 + 345t + 168)\pi_K + (750t^2 + 690t + 126))x + ((75t^2 - 120t - 30)\pi_K^4 + (-300t^2 - 120t - 15)\pi_K^3 + (-100t^2 - 29t + 46)\pi_K^2 + (1125t^2 - 1395t - 288)\pi_K + (-4500t^2 - 855t + 18))$.\\

But now $\overline{\varphi''}(x,t) = x + 3$, which is a degree 1 polynomial. Thus $\varphi(x,3+0\pi_F)= P(5,3+0\pi_F,x) = x^5 + 5x^3 - 5x^2 + 5x + 3$ has a root in $K$. And we find that $x^5 + 5x^3 - 5x^2 + 5x + 3$ is truly a Galois splitting model for $K$. \\

Similarly, for the second extension its defining polynomial is $x^5 + 10x^2 + 5$ and using Algorithm 1 we find that $t=13$ and thus its Galois splitting model is $x^5 + 10x^4 - 5x^3 + 145x^2 + 5x + 13$. For the last extension we find its defining polynomial is $x^5 + 5x^4 + 5$ and $t=18$ for a Galois splitting model of $x^5 + 15x^4 - 10x^3 + 295x^2 + 5x + 18$.\\

For a more advanced example, let $p=3$ and $K/\Q_p$ being a $C_3 \wr C_4$ extension with $F/\Q_p$ a $C_4$-subextension with a defining polynomial $f(x)=x^4 - 3x^2 + 18$. There are 16 such extensions. The field $F$ has a Galois splitting model $g(x) = x^4 + 3x^3 - 6x^2 - 18x - 9$ and $F^{unram}/\Q_p$ is a degree 2 extension with Galois splitting model $u(x) = x^2 - x - 1$. We can let a root of $\overline{u(x)}$ be the generator of $k$, and a root of $f(x)$ is the uniformizer $\pi_F$. Again from \citep{jensen2002generic} a $C_3$ generic polynomial is $Q(t,x)=x^3 - tx^2 + (t - 3)x + 1$. Using Algorithm 1 the values for the parameter of $Q(t,x)$ that will generate the Galois splitting models for all 16 $C_3 \wr C_4$ extensions where $\alpha$ is a root of $g(x)$ are listed below in Table \ref{table:table1}. \\

\begin{table}
\centering
\begin{tabular}{||p{3cm}|p{5.5cm}|p{7.5cm}||} 
 \hline
 Parameter Value & Defining Polynomial & Galois Splitting Model \\ [0.5ex] 
 \hline\hline
$-\frac{1}{3}\alpha^8- \frac{2}{3}\alpha^7+ 4\alpha^6 + \frac{10}{3}\alpha^5 - \frac{1}{3}\alpha^4 + 6\alpha^3 + 9\alpha^2$ &
$ x^{12}+ 12x^{11} - 12x^{10} - 6x^9 - 9x^7+ 6x^6+ 9x^5 + 9x^4 + 9x^3- 9$ &
$ x^{12}
 - 2673x^{11}
 + 1199940x^{10}
 - 22068644x^9
 + 91973115x^8
 - 138890646x^7
 + 55022127x^6
 + 54465408x^5
 - 53087931x^4
 + 10216039x^3
 + 1170603x^2
 + 2661x
 + 1$\\ 
 \hline
$-\frac{2}{3}\alpha^8
 - \frac{4}{3}\alpha^7
 + \frac{17}{3}\alpha^6
 + \frac{23}{3}\alpha^5
 + \frac{5}{3}\alpha^4
 + 11\alpha^3
 + 9\alpha^2$&$
 x^{12} - 6x^{10} - 3x^9 - 9x^8 + 9x^7 - 12x^6 - 9x^4 - 9$&
 $x^{12}
 - 648x^{11}
 + 53790x^{10}
 - 109544x^9
 - 251055x^8
 + 747684x^7
 - 134418x^6
 - 1037412x^5
 + 1077174x^4
 - 392936x^3
 + 46728x^2
 + 636x
 + 1$\\
 \hline
$-\frac{1}{3}\alpha^6
 - \frac{4}{3}\alpha^5
 + \frac{5}{3}\alpha^4
 + 12\alpha^3
 + 8\alpha^2$&
 $x^{12} - 27x^{11} + 21x^{10} + 39x^9 - 27x^8 + 36x^7 - 3x^6 + 36x^5- 9x^4 - 9x^3 + 27x^2 + 27x + 18$&
 $x^{12}
 - 6x^{11}
 - 6000x^{10}
 - 23405x^9
 + 148536x^8
 - 133020x^7
 - 239511x^6
 + 508500x^5
 - 332604x^4
 + 83515x^3
 - 6000x^2
 - 6x
 + 1$\\
 \hline
$\alpha^2$&$x^{12} - 9x^{11} + 12x^{10} - 9x^9+ 9x^8 - 9x^7 + 12x^6 + 9x^5 + 9x^4 + 9x^3 - 9$&$x^{12}
 - 21x^{11}
 + 135x^{10}
 - 275x^9
 - 99x^8
 + 900x^7
 - 741x^6
 - 270x^5
 + 531x^4
 - 140x^3
 - 30x^2
 + 9x
 + 1$\\
 \hline
$-\frac{1}{3}\alpha^6
 - \frac{2}{3}\alpha^5
 + 3\alpha^4
 + 5\alpha^3
 + \alpha^2$&
 $x^{12} + 24x^{11} - 39x^{10} - 3x^9 - 36x^8 + 27x^7 + 12x^6 - 18x^5+ 18x^4 + 18x^3 - 27x - 36$
 &$x^{12}
 - 48x^{11}
 - 135x^{10}
 + 886x^9
 + 90x^8
 - 3906x^7
 + 3687x^6
 + 2538x^5
 - 5436x^4
 + 2884x^3
 - 597x^2
 + 36x
 + 1$\\  
 \hline 
$-\frac{1}{3}\alpha^8
 - \frac{2}{3}\alpha^7
 + 5\alpha^6
 + \frac{14}{3}\alpha^5
 + \frac{4}{3}\alpha^4
 + 3\alpha^3
 + 5\alpha^2$&
 $x^{12} - 12x^{11} - 3x^{10} + 9x^9 + 9x^8 + 6x^6 + 9x^3 - 9$
 &$x^{12}
 - 4452x^{11}
 + 3497394x^{10}
 - 87848339x^9
 + 350436510x^8
 - 452194092x^7
 + 20067603x^6
 + 393023304x^5
 - 283549896x^4
 + 53119039x^3
 + 3448488x^2
 + 4440x
 + 1$\\
 \hline
$-\frac{1}{3}\alpha^7
 - \frac{4}{3}\alpha^6
 + \alpha^5
 + \frac{26}{3}\alpha^4
 + 14\alpha^3
 + 9\alpha^2$&$x^{12} - 6x^{11}+ 6x^{10} + 9x^9+ 9x^7 - 3x^6+ 9x^5 + 9x^4 - 9x^3 - 9$
 &$x^{12}
 - 117x^{11}
 - 5736x^{10}
 + 24646x^9
 - 12645x^8
 - 65502x^7
 + 89493x^6
 + 5544x^5
 - 67761x^4
 + 38929x^3
 - 6957x^2
 + 105x
 + 1$ \\
 \hline
$-\frac{2}{3}\alpha^7
 - \frac{4}{3}\alpha^6
 + 6\alpha^5
 + 8\alpha^4
 + 2\alpha^2$&
 $x^{12}- 3x^{11} + 6x^{10} - 12x^9 + 9x^8 - 3x^6 - 9x^3 - 9$
 &$x^{12}
 - 6x^{11}
 - 23370x^{10}
 + 70330x^9
 + 86346x^8
 - 487980x^7
 + 476484x^6
 + 70920x^5
 - 332829x^4
 + 163480x^3
 - 23370x^2
 - 6x
 + 1$\\
 \hline 
$-\frac{2}{3}\alpha^8
 - \frac{5}{3}\alpha^7
 + \frac{19}{3}\alpha^6
 + \frac{34}{3}\alpha^5
 + \frac{11}{3}\alpha^4
 + 6\alpha^3
 + 9\alpha^2$&
 $x^{12} + 3x^{11} - 6x^{10} + 12x^9 - 9x^8 - 9x^7 + 3x^6 - 9x^5 + 9x^4 - 9$
 &$x^{12}
 - 2763x^{11}
 + 1714650x^{10}
 - 24252134x^9
 + 91860525x^8
 - 127673766x^7
 + 38094477x^6
 + 61019388x^5
 - 49704831x^4
 + 7257379x^3
 + 1684323x^2
 + 2751x
 + 1$  \\ 
 \hline
$-\frac{1}{3}\alpha^7
 - \frac{4}{3}\alpha^6
 + \alpha^5
 + \frac{32}{3}\alpha^4
 + 16\alpha^3
 + 8\alpha^2$&
 $x^{12}- 30x^{11} - 33x^{10}- 3x^9 + 9x^8- 36x^7 + 30x^6+ 27x^4 + 9x^3 + 27x^2 - 27x + 18$
 &$x^{12}
 - 420x^{11}
 + 41178x^{10}
 - 140261x^9
 + 2574x^8
 + 465084x^7
 - 441147x^6
 - 239832x^5
 + 524430x^4
 - 248639x^3
 + 36624x^2
 + 408x
 + 1$  \\
 \hline
$\alpha^6
 + \frac{1}{3}\alpha^5
 + \frac{2}{3}\alpha^4
 + 6\alpha^3
 + 9\alpha^2$&
 $x^{12} + 3x^{11} + 12x^{10} - 12x^9 - 9x^8 - 12x^6 + 9x^5 - 9x^4 + 9x^3 - 9$
 &$x^{12}
 - 1998x^{11}
 + 788730x^{10}
 - 11031854x^9
 + 44262630x^8
 - 68198436x^7
 + 32763192x^6
 + 17255088x^5
 - 19860381x^4
 + 3254224x^3
 + 766818x^2
 + 1986x
 + 1$ \\
 \hline
$-\frac{1}{3}\alpha^7
 - \frac{2}{3}\alpha^6
 + \frac{7}{3}\alpha^5
 + \frac{14}{3}\alpha^4
 + 6\alpha^3
 + 8\alpha^2$&
 $x^{12} - 12x^{11} + 6x^{10} + 12x^9 - 9x^8 + 3x^6 + 9x^5 - 9x^4 - 9x^3 - 9$

 &$x^{12}
 - 366x^{11}
 + 34860x^{10}
 - 251825x^9
 + 639126x^8
 - 556380x^7
 - 147171x^6
 + 445860x^5
 - 118494x^4
 - 76865x^3
 + 30900x^2
 + 354x
 + 1$  \\
 \hline
$-\frac{2}{3}\alpha^7
 - \frac{5}{3}\alpha^6
 + \frac{14}{3}\alpha^5
 + \frac{32}{3}\alpha^4
 + 12\alpha^3
 + 10\alpha^2$&
 $x^{12} - 33x^{11} - 21x^{10} - 21x^9 - 18x^8 - 9x^7 - 24x^6 - 36x^5 + 18x^4 - 27x^3 - 27x^2 + 27x + 18$
 &$x^{12}
 - 201x^{11}
 - 22770x^{10}
 + 31720x^9
 + 347661x^8
 - 1075050x^7
 + 1010949x^6
 - 50220x^5
 - 424179x^4
 + 206815x^3
 - 24915x^2
 + 189x
 + 1$ \\ 
 \hline
$2\alpha^4
 + \alpha^2$&
 $x^{12} + 33x^{11} + 12x^{10} - 6x^9 - 18x^8 + 9x^7 + 24x^6 + 9x^5 - 36x^4 - 27x^2 + 27x - 36$
 &$x^{12}
 - 399x^{11}
 + 32373x^{10}
 - 367421x^9
 + 1045719x^8
 - 718812x^7
 - 1042023x^6
 + 1826334x^5
 - 869625x^4
 + 65416x^3
 + 28050x^2
 + 387x
 + 1$ \\
 \hline
$\alpha^4
 + \alpha^2$&
 $x^{12} + 18x^{11} + 21x^{10} - 27x^9 - 18x^8 - 12x^6 + 27x^5 - 36x^4 + 9x^3 + 27x^2 - 27x - 36$
 &$x^{12}
 - 210x^{11}
 + 9288x^{10}
 - 69071x^9
 + 146484x^8
 - 21996x^7
 - 258537x^6
 + 290628x^5
 - 91350x^4
 - 12479x^3
 + 7044x^2
 + 198x
 + 1$  \\ 
\hline
$-\frac{2}{3}\alpha^8
 - \frac{4}{3}\alpha^7
 + 7\alpha^6
 + \frac{26}{3}\alpha^5
 - \frac{2}{3}\alpha^4
 + 5\alpha^3
 + 5\alpha^2$&
 
$x^{12} + 3x^{11} + 3x^{10} - 9x^8 + 9x^7 + 3x^6 + 9x^3 - 9$
 &$x^{12}
 - 2535x^{11}
 + 1417368x^{10}
 - 8743226x^9
 + 10278819x^8
 + 21449574x^7
 - 57281997x^6
 + 41827608x^5
 - 5046435x^4
 - 5291249x^3
 + 1389549x^2
 + 2523x
 + 1$ 
\\ [1ex] 
\hline
\end{tabular}
\caption{Table showing the Galois splitting model and defining polynomial for each $C_3 \wr C_4$ extension.}
\label{table:table1}
\end{table}

\nocite{PARI2}
\nocite{bosma1997magma}

\subsection{Conclusion}
The cases in which the filtration of a Galois group by its ramification groups are interesting and most difficult and hence the need for Galois splitting models, are the wild extensions of composite degree. For wildly ramified extensions of degree 11 and lower, Jones and Roberts have computed all data relating to their ramification groups. The next interesting cases that we computed were $\mathcal{K}(2, 12)$, $\mathcal{K}(3, 12)$, $\mathcal{K}(2, 14)$, $\mathcal{K}(7, 14)$, $\mathcal{K}(3, 15)$, $\mathcal{K}(5, 15)$, and $\mathcal{K}(2, 18)$. For these cases we computed the Galois slope content, Galois mean slope, and Inertia subgroup see \citep{jj_global} on how these values are computed using the Galois splitting models. Once again all computed data are available at \url{https://hobbes.la.asu.edu/LocalFields/} and \url{http://www.lmfdb.org/LocalNumberField/} and implementation of the Algorithm 1 is located at \url{https://github.com/bcarril1/gsm_panayi}.\\

\newpage
\bibliographystyle{plain}
\bibliography{panayi-algo}

\end{document}